\documentclass{cmslatex}
\usepackage[paperwidth=7in, paperheight=10in, margin=.875in]{geometry}
\usepackage{amsmath}
\usepackage{amsfonts}
\usepackage{mathtools,pdfsync}
\usepackage{bm}

\newcommand{\fract}[2]{{#1}/{#2}}
\newcommand{\inv}{^{-1}}

\newcommand{\R}{\mathbb{R}}
\newcommand{\ds}{\,ds}

\newcommand{\e}{\varepsilon}
\newcommand{\zed}{z}
\newcommand{\imag}{\mathbf{i}}

\newcommand{\bfL}{{\bm{L}}}
\newcommand{\tfL}{{\tilde \bfL}} %{{\tilde{\kern .15em\bfL\kern -.01em}}}
\newcommand{\calI}{{\mathcal I}}

\renewcommand{\theequation}{\arabic{section}.\arabic{equation}}

	  \usepackage{hyperref}

\begin{document}
	  
          \title{Cutoff estimates for the Becker-D\"{o}ring equations}

          %For each author, make a block with the following four macros:

          \author{Ryan W. Murray\thanks{Department of Mathematics, Penn State University, State College, PA, USA. rwm22@psu.edu .} \and{ Robert L. Pego\thanks{Department of Mathematical Sciences and Center for Nonlinear Analysis, Carnegie Mellon University, Pittsburgh, PA, USA. rpego@cmu.edu .}}}

          %Use \thanks statements for acknowedgements of grants and
          %support. They will appear below all the authors' addresses, so be
          %specific about which author is thanking whom:

          \thanks{This material is based upon work supported by the National Science Foundation under grants DMS 1211161 and DMS 1515400, and partially supported by the Center for Nonlinear Analysis (CNA) under National Science Foundation PIRE Grant no.\ OISE-0967140, and the NSF Research Network Grant no.\ RNMS11-07444 (KI-Net).}

          % Use the standard latex environments for theorems, etc. Here is one
          % possible method of declaring them: It numbers all results by the
          % section, and uses a common numbering system for the different
          % environmentts.

          %\date{Received date / Revised version date}
         %\pagestyle{myheadings} \markboth{Short form for running heads}{Author's name} \maketitle
         \pagestyle{myheadings} \markboth{Cutoff for Becker-D\"oring equations}{Ryan Murray, Robert Pego}
          \maketitle

          \begin{abstract}
    This paper continues the authors' previous study (\textit{SIAM J. Math. Anal.}, 2016) of the trend toward equilibrium of the Becker-D\"{o}ring equations with subcritical mass, by characterizing certain fine properties of solutions to the linearized equation. 
In particular, we partially characterize the spectrum of the linearized operator, showing that it contains the entire imaginary axis in polynomially weighted spaces.
Moreover, we prove detailed cutoff estimates that establish upper and lower bounds 
on the lifetime of a class of perturbations to equilibrium.
%This is accomplished by partially characterizing the spectrum of the linearized operator, and by proving detailed cutoff estimates. In particular, upper and lower bounds are given on the lifetime of a class of perturbations of equilibrium.
\end{abstract}

\begin{keywords}
Coagulation-fragmentation equations, spectrum, cutoff estimates.
\end{keywords}

 \begin{AMS}
34D05. 47D06, 82C05.
\end{AMS}

          \section{Introduction.}\label{intro}
This work considers the Becker-D\"{o}ring equations, which are given by the infinite sequence of differential equations
\begin{equation} \label{OriginalBDEquation}
\begin{aligned}
\frac{d}{dt} c_i(t) &= J_{i-1}(t) -J_{i}(t),\quad i=2,3,\ldots, \\
\frac{d}{dt}  c_1(t) &= -J_1(t) - \sum_{i=1}^\infty J_i(t),
\end{aligned}
\end{equation}
where the $J_i$ can be written as
\begin{equation} \label{OriginalFluxDefinition}
J_i(t) = a_ic_1(t) c_i(t) - b_{i+1} c_{i+1}(t),
\end{equation}
and where $(a_i),(b_i)$ are fixed, positive sequences, known as the coagulation and fragmentation coefficients respectively. 
These equations are a well-known model for certain physical phenomena 
occuring in phase transitions, such as condensation in alloys and polymers. 
In this context, $c_i(t)$ typically represents the density of particles of size $i$ 
(``$i$-particles'') in some units.
The Becker-D\"oring equations \eqref{OriginalBDEquation} describe the evolution of the discrete size distribution $(c_i)$ under mean-field assumptions 
which state that $i$-particles aggregate with $1$-particles (monomers)
to form $i+1$-particles at rate $a_i c_1(t)$ per particle, 
and $i+1$-particles break in two pieces, monomers and $i$-particles, 
at rate $b_{i+1}$ per particle.
The first moment $\mu=\sum_{i=1}^\infty ic_i(t)$ 
corresponds to the total mass in the system and  
is formally conserved in time, due to the evolution equation for $c_1(t)$ in \eqref{OriginalBDEquation}. 

The quantity $J_i(t)$ is the net reaction rate for $i$-particles to become $i+1$-particles,
and this vanishes in equilibrium. Under typical assumptions on the rate coefficients 
(described below), it is known \cite{BallCarrPenrose} that 
there is a critical mass $\mu_{\rm crit}\le\infty$ such
that for positive initial data $(c_i^0)$ with subcritical mass, meaning
\begin{equation}\label{eqn:subcritical-data}
\sum_{i=1}^\infty i c_i^0 =: \mu < \mu_{\rm crit}\,,
\end{equation}
the solution $(c_i(t))$ converges strongly to an equilibrium solution $(Q_i)$ 
determined by the condition that $J_i=0$ for all $i$, i.e.,
\begin{equation}\label{DetailedBalance}
Q_1 = z, \qquad b_{i+1}Q_{i+1} = a_i Q_i Q_1, \qquad \sum_{i=1}^\infty Q_i i = \mu.
\end{equation}
Here by strong convergence we mean that
\begin{equation} \label{eqn:strong-convergence}
\lim_{t \to \infty} \sum_{i=1}^\infty i|c_i(t) - Q_i| =0.
\end{equation}

Various authors have sought to establish uniform convergence rates in \eqref{eqn:strong-convergence}. The previous works \cite{CanizoLods,JabinNiethammer} focused on convergence rates in the setting where the initial data decays exponentially fast. More recent works \cite{CanizoLodsEinav,MurrayPego2015} have focused on convergence rates when the initial data decays only algebraically fast. In particular, in \cite{MurrayPego2015} the present authors proved the following result.

\begin{theorem}[\cite{MurrayPego2015}] \label{thm:1}
Assume the model coefficients in \eqref{OriginalFluxDefinition} satisfy conditions \eqref{aLowerBound}-\eqref{boundedByI} below. 
Let $(c_i(t))$ be a solution of the Becker-D\"oring equations \eqref{OriginalBDEquation} 
with subcritical mass, and let its deviation from equilibrium be represented by
$(h_i(t))$, defined so that
\begin{equation}\label{hDefinition}
c_i = Q_i (1 + h_i).
\end{equation}
Let  $m$ and $k$ be real numbers satisfying $m> 0$ and $k>m+2$. Then there exists positive constants $\delta_{k,m}, C_{k,m}$ so that if $\|h(0)\|_{X_{1+k}} < \delta_{k,m}$ then 
\begin{equation}\label{eqn:Alg-Decay}
\|h(t)\|_{X_{1+m}} \leq C_{k,m}(1+t)^{-(k-m-1)}\|h(0)\|_{X_{1+k}} \quad \mbox{for all $t \geq 0$}.
\end{equation}
\end{theorem}

Here we are writing
\begin{equation}
X_k := \left\{ (h_i) : \|h\|_{\ell^1(Q_i i^k)} := \sum_{i=1}^\infty Q_i i^k |h_i| < \infty,\quad  \sum_{i=1}^\infty Q_i i h_i = 0\right\},\quad k\geq 1,\\
\end{equation}
with norm $\|\cdot\|_{X_k}=\|\cdot\|_{\ell^1(Q_i i^k)}$.

Theorem~\ref{thm:1} was derived from analysis conducted on the linearized equation
\begin{equation}\label{eqn:linear}
\frac{d}{dt} h = Lh, 
\end{equation}
where the operator $L$ is defined in weak form on suitable spaces by the requirement that for all suitable
test sequences $(\phi_i)$,
\begin{equation}\label{d:Lweak}
\sum_{i=1}^\infty Q_i (Lh)_i \phi_i = \sum_{i=1}^\infty Q_i Q_1 a_i (h_{i+1}-h_i-h_1)(\phi_1 + \phi_i - \phi_{i+1}).
\end{equation}
In particular, it was first shown that $e^{Lt}$ is uniformly bounded in $X_1$, after which bounds of the type \eqref{eqn:Alg-Decay} were obtained for the linearized equation via interpolation theory.
%, see \cite{MurrayPego2015} for details. 
The linearized estimates were then extended, after some technicalities, to the non-linear setting. 

A natural question is whether the bounds in this theorem are optimal. More generally, one would hope for a more detailed understanding of the dynamics as solutions converge to equilibrium. 
The first theorem in the present work seeks to address these questions by giving detailed information about the spectrum of the linearized operator.

\begin{theorem}\label{thm:spectrum}
Suppose, in addition to \eqref{aLowerBound}-\eqref{boundedByI}, that 
\[ \mbox{$a_i - a_{i-1} = o(1)$, \quad $b_i - b_{i-1} = o(1)$ \quad and \quad $a_i \to \infty$. }\]
Then the operator $L$ has an approximate eigenvalue at $\lambda \imag$ in $X_k$ for all $\lambda \in \R$ and all $k \geq 1$, where $\imag = \sqrt{-1}$.
\end{theorem}

Here we write $\imag= \sqrt{-1}$ so that we may use $i$ freely as an index throughout the work. We remark that the assumptions of this theorem are satisfied by a wide class of coefficients used in applications, see e.g. \eqref{PenroseCoefficients}.

Theorem~\ref{thm:spectrum} highlights significant differences between the operator in exponentially weighted spaces as opposed to polynomially weighted ones. In exponentially weighted spaces one finds that the operator $L$ generates an analytic semigroup with uniform decay, and can even be self-adjoint \cite{CanizoLods}. On the other hand, in polynomially weighted spaces, the operator $L$ only generates a bounded $C_0$ semigroup.
The hypotheses of Theorem~\ref{thm:spectrum} cover most physically relevant
cases, for which $a_i \sim i^\alpha$ with $\alpha \in (0,1)$. In these cases,
Theorem~\ref{thm:spectrum} shows the operator $L$ actually has approximate
point spectrum at every point on the imaginary axis, thus it cannot generate an
analytic semigroup.
%This demonstrates a significant departure from the analytic semigroup setting.

It is our viewpoint that this spectral phenomenon may give insight into some of the difficulties encountered in the study of coagulation-fragmentation equations in a more general setting. Here, in the most natural space $\ell^1(i)$, the operator $L$ has spectrum on the imaginary axis, suggesting that systems of this type should be treated like hyperbolic equations. Indeed, many of the techniques used in studying coagulation-fragmentation equations (such as entropy methods and limits of regularizations \cite{BallCarrPenrose}) originate in the study of hyperbolic equations. We believe that some of the obstacles for analyzing coagulation-fragmentation dynamics may directly relate to ``hyperbolic'' aspects of the equations.

The next theorem seeks to give more detailed information on the dynamics of the linearized system 
with `pulse-like' initial data supported far from the origin, 
with $0\ll N_1<i<N_2$.  Below, the \textit{support} of $h^0$ is the set $\{i:h^0_i\neq 0\}$.

%\begin{theorem}\label{thm:dynamics}
%Suppose, in addition to \eqref{aLowerBound}-\eqref{boundedByI}, that $a_i = i^\alpha$ with $\alpha \in (0,1)$. Suppose that $h^0$ is a positive sequence, with support between values $N_1$ and $N_2$, satisfying $\sum Q_i i h_i = 1$. Fix $\e >0$. Then there exists a $K$ and $\delta$ so that $h(t) = e^{Lt}h^0$, and for times $t < T = \delta N_1^{1-\alpha}$ we have that $\sum Q_i i h_i(t) \chi(i,t) > 1-\e$, where $\chi(t)$ has the form
%\[
%\chi(i,t) = \begin{cases} 1 &\text{ if } A(2t,N_1)-K < i < A(t/2,N_2) + K,\\ 0 &\text{ otherwise}, \end{cases}
%\]
%where $A_1,A_2$ are solutions of a specific ODE, see \eqref{eqn:A-def}.
%\end{theorem}

\begin{theorem}\label{thm:dynamics}
Suppose, in addition to \eqref{aLowerBound}-\eqref{boundedByI}, that ${a_i}/{i^\alpha} \to 1$ with $\alpha \in (0,1)$. Let $h^0=(h_i^0)$ be a non-negative sequence satisfying $\sum Q_i i h_i^0 = 1$ and 
let 
\[
h(t) = e^{Lt} h^0
\]
be the solution of \eqref{eqn:linear} in $\ell^1(Q_i i)$ with initial data $h^0$.
For any $\e>0$, there exist $N^*$, $K^*$ and $\delta>0$
%Fix $\e > 0$. Then there exists a $K_1$, $K_2$ and $\delta$ 
    such that whenever $h^0$ is supported in $\{N_1<i<N_2\}$ with
    $N^*<N_1<N_2$, 
   % unless $N_1<i<N_2$ (with $K_1 < N_1 < N_2$), 
    then for all times $t < T = \delta N_1^{1-\alpha}$ we have that 
\[
\sum_{i=1}^\infty Q_i i h_i(t) \chi(i,t) > 1-\e,
\]
 where $\chi$ has the form
\[
\chi(x,t) = \begin{cases} 1 &\text{ if } A(N_1,2t)-K^* < x < A(N_2,t/2) + K^*,\\ 0 &\text{ otherwise}, \end{cases}
\]
where $A(x,t)$ is the solution of 
\begin{equation}\label{eqn:A-def}
\frac{\partial}{\partial t} A = -(z_s-z) A^\alpha, \quad A(x,0) = x.
\end{equation}
Furthermore $N^*$ can be chosen independent of $\e, \delta, $ and $K^*$.
%Furthermore, $K_1$ can be chosen independent of $\e, N_1$.
\end{theorem}

In \eqref{eqn:A-def}, $z_s$ is the critical monomer density, defined in \eqref{tildeQLimit} below.
The explicit solution of \eqref{eqn:A-def} is 
\begin{equation}\label{eqn:A-sol}
A(x,t) = (x^{1-\alpha} - (z_s-z) (1-\alpha)t)^{1/(1-\alpha)}.
\end{equation}
%
%\begin{equation}\label{eqn:A-def}
%\begin{aligned}
%\frac{d}{dt} A_1 &= -2 c_2 A_1^\alpha, \quad  A_1(0,x) = x, \\
%\frac{d}{dt} A_2 &= -\frac{c_1 A_2^\alpha}{2}, \quad  A_2(0,x) = x,
%\end{aligned}
%\end{equation}
%where $c_1$ and $c_2$ are given by
%\[
%c_1i^\alpha \leq  \left(\frac{a_i Q_1 Q_i i}{(i+1)Q_{i+1}} - \frac{b_i i Q_i}{(i-1)Q_{i-1}}\right) \leq c_2 i^\alpha.
%\]
The intuition behind the result in Theorem~\ref{thm:dynamics} can be explained as follows. After writing $u_i = Q_i i h_i$, we can formally approximate $e^{Lt}u$ by solving an advection diffusion equation of the form
\[
u_t = p(x) u_{x} + q(x) u_{xx},
\]
where $p(x)$ grows like $(z_s - z) x^\alpha$ and $q(x)$ grows like $z x^\alpha$. If we neglect the diffusion term, we find that the solution is constant along the characteristic curves precisely given by $A$, and $\chi$ then describes how ``mass'' travels through the system. Thus the result of Theorem~\ref{thm:dynamics} essentially tells one that this advection is sufficient to describe the spread of the mass, at least when pulses are well-separated from $i = 1$.

Of course there are error terms in approximating the evolution of $(u_i)$ by this advection-diffusion equation, but these error terms go to zero for large-enough cluster size $i$. We also remark that we do not in fact solve the advection-diffusion equation in proving Theorem~\ref{thm:dynamics}; instead we opt to construct a supersolution and then work exclusively with the discrete equations.

One natural application of the previous theorem is the following corollary.

\begin{corollary}\label{cor:cutoff}
Suppose, in addition to \eqref{aLowerBound}-\eqref{boundedByI}, that $a_i / i^\alpha\to1$ with $\alpha \in (0,1)$. Let $N^*$ be given as in the statement of Theorem \ref{thm:dynamics}. Then for any $\e > 0$ there exists  $\delta >0$ such that for any $N > N^*$ there exists an $(h_i^0) \in X_1$ with support in $\{i < N\}$ and with $\|h^0\|_{X_1} = 1$ satisfying
\begin{equation}\label{eqn:lower-cutoff}
\|e^{Lt} h^0\|_{X_1} \geq 1-\e
\end{equation}
for all $t < \delta N^{1-\alpha}$. In particular, $\|e^{Lt}\|_{\mathcal{L}(X_1)} \geq 1$ for all $t \geq 0$.

On the other hand, for any $\e > 0$, there exists a $\delta > 0$ so that for any $(h_i^0) \in X_1$ with support in $\{i < N\}$ satisfying $\|h^0\|_{X_1} = 1$ we have that
\begin{equation} \label{eqn:upper-cutoff}
\|e^{Lt} h^0\|_{X_1} \leq \e
\end{equation}
 for all $t > \delta N$.
\end{corollary}

This corollary implies that the estimates in \cite{MurrayPego2015} are optimal in the sense that we cannot expect any uniform decay estimates for data in the natural space $X_1$. We remark that the upper bound \eqref{eqn:upper-cutoff} on perturbation lifetimes is a direct consequence of the decay estimates in exponentially weighted spaces derived in \cite{CanizoLods}, and may not be sharp. The novel contribution here is the lower bound \eqref{eqn:lower-cutoff}.

The result of Corollary~\ref{cor:cutoff} is analogous to the \emph{cutoff phenomenon} in the theory of Markov chains. In short, a Markov chain is said to exhibit a cutoff phenomenon if typical states remain far from equilibrium up to well-quantified time after which a rapid transition to equilibrium occurs.
See \cite{DiaconisReview} for a detailed introduction to the subject. Examples of Markov chains exhibiting such behavior include card shuffling\footnote{This is easily remembered by the rule of thumb given in \cite{BayerDiaconis} that a deck of cards is not random until shuffled 7 times.}\cite{BayerDiaconis} and random walks on a hypercube \cite{DiaconisGrahamMorrison}. If we view the linearized Becker-D\"oring equations as a continuous time Markov chain, then the inequalities in Corollary \ref{cor:cutoff} precisely describe a type of cutoff phenomenon, where the cutoff time depends on the support of the initial data.

We remark that in studying cutoff phenomena, the norm that is used is often critically important. For example, in the case of random walks on a hypercube, if deviations from equilibrium are measured in $\ell^2$ then there is no cutoff (in fact the transition matrix is symmetric), but if measured in $\ell^1$ a cutoff phenomenon occurs \cite{JonssonTrefethen, TrefethenBook}. The analogy continues to hold for the linearized Becker-D\"oring equations: $e^{Lt}$ is exponentially decaying (and $L$ is self-adjoint) in exponentially-weighted $\ell^2$ spaces but displays a persistence phenomenon in polynomially-weighted $\ell^1$ spaces.

The analogies actually go even deeper. In \cite{JonssonTrefethen}, Jonsson and Trefethen suggest that for the random walk on a hypercube the 1-pseudospectra help explain the cutoff phenomenon. In the Becker-D\"oring case, Theorem \ref{thm:spectrum} establishes the existence of approximate point spectrum (in the 1-norm) on the imaginary axis. Jonsson and Trefethen also explain the cutoff phenomenon for the random walk on a hypercube case in terms of the overlap of two sliding Gaussians, or in other words in terms of an advection phenomenon. Our results above show that the cutoff times in the Becker-D\"oring case are similarly explained by advection. In short, the linearized Becker-D\"oring equations exhibit many of the same features seen in Markov chains that exhibit cutoff.

At this point, we are careful to remark that our cutoff result is not sharp, in the sense that the upper and lower bounds do not match. More delicate analysis would be required to obtain sharp results in this direction.

We also note that our results in this work do not address the cutoff phenomenon in
the context of nonlinear Becker-D\"oring dynamics.
Our reason for focusing on the linear case is that it highlights
some fundamental obstacles to uniform rates of convergence and the
differences between the dynamics in different function spaces, without
getting too bogged down in technicalities.

Consistent with other works on the Becker-D\"oring equations \cite{JabinNiethammer}, \cite{MurrayPego2015}, we will make the following standard assumptions on the model coefficients $(a_i), (b_i)$:
\begin{align}
a_i >C_1&>0\qquad \text{ for all } i \geq 1, \label{aLowerBound}\\
\lim_{i \to \infty} \frac{a_{i+1}}{a_i} &= 1, \label{aLimit}\\
\lim_{i \to \infty} \frac{a_i}{b_{i}} &=: \frac{1}{\zed_s}\in(0,\infty) \label{tildeQLimit}\\
 a_i, b_i &\leq C_2 i   \qquad \text{ for all } i\geq 1. \label{boundedByI}
\end{align}

These assumptions are satisfied by many of the coefficients proposed for physical phenomenon.
For example, these assumptions are satisfied by the coefficients proposed in \cite{Penrose89}
\begin{equation}\label{PenroseCoefficients}
a_i = i^\alpha,\quad b_i = a_i\left(z_s + \frac{q}{i^{1-\beta}}\right),\quad \alpha \in (0,1],\quad \beta \in [0,1],\quad q > 0.
\end{equation}

Using  \eqref{DetailedBalance}, \eqref{aLimit}, and \eqref{tildeQLimit}, it is straightforward to show that
\begin{equation}
\lim_{i \to \infty} \frac{Q_{i+1}}{Q_i} = \frac{z}{z_s}. \label{QLimit}
\end{equation}
In fact, if we write $Q_i = \tilde Q_i z^i$, then $z_s$ and $\mu_s$ are given by
\[
z_s = \sup \left\{z : \sum_{i=1}^\infty i \tilde Q_i z^i < \infty \right\}, \quad \mu_s = \sup\left\{\sum_{i=1}^\infty i \tilde Q_i z^i : z < z_s \right\}. 
\]
Hence the restriction to subcritical data, that is \eqref{eqn:subcritical-data}, implies that $\frac{z}{z_s} < 1$, which in turn implies that the $Q_i$ are exponentially decaying.

\section{Approximate Spectrum.}
The aim of this section is to prove Theorem \ref{thm:spectrum}. The main idea is to construct approximate eigenvectors using wide pulses of constant modulus with support far away from $i=1$. A simple version of the theorem, in the case $\lambda = 0$, can be found in the first author's thesis, see Theorem 8.2.10 in \cite{MurrayThesis}.

%\begin{proof}[Proof of Theorem \ref{thm:spectrum}]
\textit{Proof of Theorem \ref{thm:spectrum}.}
Define $\tilde h_i$ so that
\[
iQ_i\, \tilde h_i = \begin{cases} 0 &\text{ if } i < N_1\\
{\exp \left( \lambda \imag(z_s-z)^{-1} \sum_{j=N_1}^i a_j^{-1} \right) } &\text{ if } N_1 \leq i \leq N_2 \\ 
% \frac{\exp \left( \lambda \imag(z_s-z)^{-1} \sum_{j=N_1}^i a_j^{-1} \right) }{i Q_i} &\text{ if } N_1 \leq i \leq N_2 \\ 
 0 &\text{ if } N_2 < i \end{cases}
\]
where $N_1 < N_2$ are constants to be determined. Clearly
\begin{equation} \label{eqn:total-mass-1}
\sum_{i=1}^\infty Q_i i^k |\tilde h_i| = \sum_{i=N_1}^{N_2} i^{k-1}.
\end{equation}
Furthermore, for $N_1 < i < N_2$, letting $w_i := \exp \left( \lambda \imag(z_s-z)^{-1} \sum_{j=N_1}^i a_j^{-1} \right)$,
\begin{align*}
&Q_i i^k ((L \tilde h)_i-\lambda \imag\tilde h_i) = i^k Q_i \left( b_i (\tilde h_{i-1} - \tilde h_i) + a_i Q_1 (\tilde h_{i+1} - \tilde h_i)- \lambda \imag\tilde h_i \right) \\
&= i^{k-1}  w_i  \left(b_i\left( \frac{Q_i i\exp \left(-\frac{\lambda \imag}{(z_s-z)a_i} \right)}{Q_{i-1}(i-1)} - 1\right) + a_i Q_1 \left(\frac{Q_i i\exp \left(\frac{\lambda \imag}{(z_s-z) a_{i+1}}\right)}{Q_{i+1} (i+1)} -1 \right) - \lambda \imag\right) \\
&= i^{k-1} w_i \left(-b_i + a_{i-1}Q_1 \frac{ i\exp \left(-\frac{\lambda \imag}{(z_s-z) a_i} \right)}{(i-1)} - a_i Q_1 + b_{i+1}\frac{ i\exp \left(\frac{\lambda \imag}{(z_s-z) a_{i+1}}\right)}{ (i+1)}  -\lambda \imag\right),
\end{align*}
where we have used \eqref{DetailedBalance}. By using the Taylor expansion of $\exp(w)$ near $w = 1$, and recalling that $\frac{i}{i+1} = 1 + O(i^{-1})$, we may use \eqref{tildeQLimit} to find that
\begin{align*}
&Q_i i^k ((L \tilde h)_i-\lambda \imag\tilde h_i) \\
&\quad = i^{k-1}  w_i \left( a_{i-1}Q_1 - b_i + b_{i+1} -a_iQ_1 + \lambda \imag\left(\frac{-a_{i-1}Q_1}{a_i(z_s-z)} + \frac{b_{i+1}}{a_{i+1}(z_s-z)} - 1 \right) 
\right.    \\
    &\quad  \left.\qquad +\  O\left(\frac{1}{a_i}\right) + O\left(\frac{a_i}{i}\right)  \right).
\end{align*}
%
%&= i^{k-1} w_i \left(-b_i + a_{i-1}Q_1 \exp \left(-\frac{\lambda \imag(z_s-z)^{-1}}{a_i} \right) - a_i Q_1 + b_{i+1}\exp \left(\frac{\lambda \imag(z_s-z)^{-1}}{a_{i+1}}\right) -\lambda \imag+ O\left( \frac{a_i}{i}\right) \right) \\
%&= i^{k-1}  w_i \left( a_{i-1}Q_1 - b_i + b_{i+1} -a_iQ_1 + \lambda \imag\left(\frac{-a_{i-1}Q_1}{a_i(z_s-z)} + \frac{b_{i+1}}{a_{i+1}(z_s-z)} - 1 \right) +  O(a_i^{-1}) + O(\frac{a_i}{i})  \right).
%\end{align*}
Since $a_i-a_{i-1} = o(1)$, we remark that
\begin{equation} \label{eqn:sublinear}
\frac{a_i}{i} = \frac{a_1 + \sum_{j=1}^{i-1}(a_{j+1}-a_{j})}{i} \to 0.
\end{equation}
Hence, using the assumptions that $a_i-a_{i-1} = o(1)$, $b_i - b_{i-1} = o(1)$ and $a_i \to \infty$ we find that
\[
Q_i i^k ((L \tilde h)_i-\lambda \imag\tilde h_i) =  i^{k-1}  w_i \left(\lambda \imag\left(\frac{-a_{i-1}Q_1}{a_i(z_s-z)} + \frac{b_{i+1}}{a_{i+1}(z_s-z)} - 1 \right) + o(1) \right)
\]
Recalling \eqref{aLimit}, \eqref{tildeQLimit}, and that $Q_1 = z$, we find that for any $\delta > 0$ we may choose $N_1$ large enough that
\begin{equation}\label{eqn:est1}
Q_i i^k |(L \tilde h)_i-\lambda \imag\tilde h_i| \leq \delta i^{k-1}
\end{equation}
for all $N_1 < i < N_2$.

On the other hand, for any $i > 1$ we have, by \eqref{tildeQLimit} and \eqref{QLimit}, that
\begin{align*}
Q_i i^k |(L \tilde h)_i-\lambda \imag\tilde h_i |  &=  Q_i i^k \left| b_i (\tilde h_{i-1} - \tilde h_i) + a_i Q_1 (\tilde h_{i+1} - \tilde h_i) -\lambda \imag\tilde h_i  \right| \\
&\leq i^k \left( \left|\frac{b_i Q_i}{Q_{i-1}(i-1)}\right| + \left|\frac{b_i}{i}\right| + \left|\frac{a_iQ_1Q_{i}}{Q_{i+1}(i+1)}\right| + \left|\frac{a_iQ_1}{i}\right|  + \left|\frac{\lambda}{i}\right| \right)\leq  C i^{k-1}a_i,
\end{align*}
where $C$ is independent of $i, N_1,$ and $N_2$. Equation \eqref{eqn:sublinear} then implies that for any $\delta>0$ we may choose a $N_1$ large enough that for all $i>1$ we have
\begin{equation}\label{eqn:est2}
Q_i i^k |(L \tilde h)_i-\lambda \imag\tilde h_i | < \delta i^k.
\end{equation}
Similarly, for $i=1$ by \eqref{aLimit} and \eqref{QLimit} we find that
\begin{align*}
|Q_1 (L \tilde h)_1 - \lambda \imag\tilde h_1| &= \left|\sum_{j=1}^\infty a_j Q_j Q_1 (\tilde h_{j+1}-\tilde h_j)\right| \\
&\leq C \sum_{j=N_1}^{N_2} \frac{a_j}{j} ,
\end{align*}
with $C$ independent of $N_1,N_2$. Hence, by \eqref{eqn:sublinear}, for any $\delta > 0$ we may again choose $N_1$ large enough that
\begin{equation}\label{eqn:est3}
|Q_1 (L \tilde h)_1 - \lambda \imag\tilde h_1| < \delta (N_2-N_1).
\end{equation}
Choosing $N_2 = 2N_1$ it is straightforward to show that
\begin{equation} \label{eqn:total-mass}
C_1 N_1^k < \sum_{i=1}^\infty i^{k-1}  \leq C_2 N_1^k.
\end{equation}
On the other hand, by equations \eqref{eqn:est1},\eqref{eqn:est2}, \eqref{eqn:est3} and \eqref{eqn:total-mass}, we have that
\begin{align*}
\sum_{i=1}^\infty Q_i i^k |(L \tilde h)_i - \lambda \imag\tilde h_i| &\leq |Q_1 (L \tilde h)_1 - \lambda \imag\tilde h_1| + |Q_{N_1-1} (L \tilde h)_{N_1-1} - \lambda \imag\tilde h_{N_1-1}| \\
&+ |Q_{N_1} (L \tilde h)_{N_1} - \lambda \imag\tilde h_{N_1}|  + |Q_{N_2} (L \tilde h)_{N_2} - \lambda \imag\tilde h_{N_2}| \\
&+ |Q_{N_2+1} (L \tilde h)_{N_2+1} - \lambda \imag\tilde h_{N_2+1}| + \sum_{i=N_1+1}^{N_2-1} Q_i i^k |(L \tilde h)_i - \lambda \imag\tilde h_i| \\
&\leq  \delta N_1 + 4 \delta N_1^k + \delta C_2 N_1^k,
\end{align*}
where $\delta \to 0$ as $N_1 \to \infty$. On the other hand, by \eqref{eqn:total-mass-1} and \eqref{eqn:total-mass} we have that
\[
\sum_{i=1}^\infty Q_i i^k |\tilde h_i| \geq C_1 N_1^k.
\]
Hence we find that
\[
\lim_{N_1 \to \infty} \frac{\sum_{i=1}^\infty Q_i i^k |(L \tilde h)_i - \lambda \imag\tilde h_i |}{\sum_{i=1}^\infty Q_i i^k |\tilde h_i|} = 0.
\]
By summing two of these pulses with non-overlapping support, scaled so that the mass constraint is satisfied, we obtain the desired result. This completes the proof.
\endproof
%\end{proof}

\section{Cutoff Phenomenon.}
The aim of this section is to prove Theorem \ref{thm:dynamics}. Before we begin the proof, we will give some definitions and recall key facts.

We recall that $L$ generates a semigroup of contractions on $\ell^2(Q_i)$ (see Section 2 in \cite{CanizoLods}). We also recall that $L$ generates a bounded semigroup on $X_1$, namely the zero mass subspace of $\ell^1(Q_i i)$ (see Theorem 2.11 in \cite{MurrayPego2015}). 
In the proof, it will be necessary to consider semigroups on the space $\ell^1(Q_ii)$, not $X_1$. To this end, note that $\xi_i = {i}/{\sum_{i=1}^\infty Q_i i^2}$ is an eigenvector (with eigenvalue $0$) of the operator $L$, normalized in $\ell^1(Q_i i)$. Furthermore, the linear mapping 
\[
h\mapsto \mu(h) = \sum_{i=1}^\infty Q_i i h_i
\] is continuous on the space $\ell^1(Q_i i)$. Hence for $h \in \ell^1(Q_i i)$ we can write
\[
e^{Lt} h = e^{Lt} (h - \xi \mu(h)) + \xi \mu(h).
\]
Because $h - \xi \mu(h)\in X_1$, 
it is then straightforward to estimate
\begin{align*}
\|e^{Lt} h\|_{\ell^1(Q_i i)} &\leq \|e^{Lt} (h - \xi \mu(h))\|_{\ell^1(Q_i i)} + \|\xi \mu(h)\|_{\ell^1(Q_i i)} \\
&\leq M\|(h - \xi \mu(h))\|_{\ell^1(Q_i i)} + C\|h\|_{\ell^1(Q_i i)}  \leq C\|h\|_{\ell^1(Q_i i)} .
\end{align*}
We also remark that $e^{Lt}$ preserves mass, in the sense that for any $h \in \ell^1(Q_i i)$ 
\begin{equation}\label{eqn:mass-preserved}
\sum_{i=1}^\infty Q_i i (e^{Lt} h)_i = \sum_{i=1}^\infty Q_i i h_i.
\end{equation}

The strategy in proving Theorem \ref{thm:dynamics} is to approximate the operator $L$ by a tridiagonal operator $\tilde L$ which is in ``divergence'' form
in mass-weighted variables. 
From \eqref{d:Lweak}, we note that for $i > 1$,
\[
iQ_i ( L h)_i = iQ_i Q_1 a_i (h_{i+1} - h_i - h_1) - iQ_{i-1} Q_1 a_{i-1} (h_i-h_{i-1}-h_1),
\]
while for $i =1$,
\[
Q_1 (Lh)_1 = Q_1^2a_1(h_2-2h_1)+\sum_{i=1}^\infty Q_i Q_1 a_i (h_{i+1} - h_i - h_1).
\]
We define the operator $\tilde L$ by requiring that for $i>1$,
\[
iQ_i  (\tilde L h)_i = 
 iQ_i Q_1 a_i (h_{i+1} - h_i) - (i-1)Q_{i-1} Q_1 a_{i-1} (h_i-h_{i-1}),
%Q_{i-1} Q_1 a_{i-1} \frac{i-1}{i} h_{i-1} - \left(Q_{i-1}Q_1 a_{i-1} \frac{i-1}{i} + Q_iQ_1a_i\right) h_i + Q_i Q_1 a_{i} h_{i+1},
\]
and for $i=1$,
\[
Q_1 (\tilde L h)_1 =  Q_1^2 a_{1}(h_2-h_1).
%Q_1 (\tilde L h)_1 := -Q_1^2a_1 h_1+ Q_1^2 a_{1} h_2.
\]

We remark that $L -\tilde L$ is a bounded operator in $\ell^2(Q_i)$, since
\begin{align*}
\|(L - \tilde L) h\|_{\ell^2(Q_i)} &\leq \left( \sum_{i=1}^\infty Q_i \left( (b_i - Q_1a_i)h_1 + \frac{b_i}{i} h_{i-1} - \frac{b_i}{i} h_i\right)^2\right)^{1/2} \\
&+ Q_1^{1/2}\left| \sum_{i=1}^\infty Q_i a_i (h_{i+1} - h_i - h_1) \right| + Q_1^{1/2} (Q_1 a_1(h_2 - 2h_1))^2  \\
&\leq C \|h\|_{\ell^2(Q_i)}^2,
\end{align*}
where we have used the assumption that $\fract{b_i}{i} \to 0$, the Cauchy-Schwarz inequality, and the fact that $\sum_{i=1}^\infty Q_i a_i^2 < \infty$.

It will be more convenient throughout the proof to work with the mass-weighted variables
defined by  $v_i = Q_i i h_i$, so that
\[
\sum_{i=1}^\infty Q_i i |h_i| = \sum_{i=1}^\infty |v_i|.
\]
We can write $v=\mathcal I h$ in terms of the diagonal operator $\mathcal I$ with entries $Q_i i$. Then we can express the operator $L$ in these new coordinates via similarity transformation as $\bfL v := \mathcal I L \mathcal I^{-1} v$. Explicitly, 
the operator $\bfL$ is given by 
\begin{equation}\label{eqn:advection-term}
(\bfL v)_i = i(a_{i-1}Q_{i-1} - a_i Q_i)v_1 + \left( -a_i Q_1-b_i \right) v_i  + a_{i-1}Q_1\frac{ i}{i-1}v_{i-1}  + b_{i+1} \frac{i}{i+1} v_{i+1} ,
\end{equation}
for $i>2$, and for $i=1$ via
\[
(\bfL v)_1 =  \sum_{i=1}^\infty \frac{b_{i+1}}{i+1} v_{i+1} - \frac{a_i Q_1}{i}v_i - a_i Q_i v_1.
\]
Similarly, in these coordinates, letting $\tfL v := \mathcal I \tilde L \mathcal I^{-1} v$,  we find that
\begin{equation} \label{eqn:tilde-L}
\begin{aligned}
(\tfL v)_i &= a_{i-1}Q_1 v_{i-1} - a_{i} Q_1 v_{i}  
+  b_{i+1}  \frac{i}{i+1}v_{i+1} - b_i \frac{i-1}{i}v_{i} , \quad \text{ for } i > 1, \\
(\tilde\bfL v)_{1} &= -a_{1} Q_1 v_1 + b_{2} \frac{1}{2} v_{2}.
\end{aligned}
\end{equation}
These expressions show that $\tilde\bfL v$ takes the form of a discrete
``divergence,'' a fact that will be useful below. 

Given initial data $u^0$ with fixed support, we will show how mass is advected through the system using a type of minimum principle. In particular, we prove a minimum principle on the ``integrated'' operator, that is the partial sums of $e^{\tfL t} u^0$. After obtaining good controls on the support of $e^{\tfL t} u^0$, we then use Duhamel's formula to establish suitable estimates on $(e^{\bfL t}- e^{\tfL t})u^0$.

%\begin{proof}[Proof of Theorem \ref{thm:dynamics}]
\textit{Proof of Theorem \ref{thm:dynamics}}.
%We note that
%\[
%\partial_t i c_i = J_{i-1} + (i-1) J_{i-1} - i J_i.
%\]
%Our goal will be to establish suitable estimates on $u_i = Q_i i h_i$, which are equivalent to $\ell^1$ estimates. We will accomplish this by proving some bounds on the solution to the equation
%\[
%\partial_t v_i = (i-1) J_{i-1} - i J_i.
%\]
%If we sum this equation, we obtain
%\[
%\partial_t \sum_{i=1}^N v_i \approx N J_N.
%\]
Given initial data $(h^0)$ as in the statement of the theorem, we set 
\[
h = e^{Lt} h^0,\qquad u^0 = \mathcal I h^0 =  (Q_i i h_i^0),\qquad 
u(t) = \mathcal I h(t)=e^{\bfL t} u^0.
\]
Clearly
\[
\sum_{i=1}^\infty  |u_i| = \sum_{i=1}^\infty Q_i i |h_i|
\]
Given initial data $u^0$ as in the assumptions, we let 
\[
v(t) := e^{\tfL t} u^0 = \calI e^{\tilde L t} \calI\inv u^0,
\]
 and set 
 \[
 V_i := \sum_{j=1}^i v_i.
 \]
  We find that the following identity holds, for all $i\ge1$:
\begin{equation}\label{eqn:integrated-form}
\begin{aligned}
\frac{d}{dt} V_i &= -a_iQ_1 v_i(t) + b_{i+1} v_{i+1}(t) \frac{i}{i+1}\\
& =   -a_iQ_1 (V_i - V_{i-1}) + b_{i+1} (V_{i+1} - V_i) \frac{i}{i+1} =: (\mathbb{L}V)_i,
\end{aligned}
\end{equation}
where we let $V_0 = 0$, and we have used the fact that $V_i - V_{i-1} = v_i$.

Next, we claim that if $u^0$ has compact support then $\lim_{i\to \infty} V_i(t) = 1$ for all $t$. As $L$ generates a semigroup of contractions on $\ell^2(Q_i)$, and as $\tilde L$ is a bounded perturbation of $L$ in $\ell^2(Q_i)$, we then have that $e^{\tilde L t}$ also generates a semigroup on $\ell^2(Q_i)$ with bound $Ce^{Ct}$. Thus for any $t>0$, $e^{\tilde L t}h^0$ is an element of $\ell^2(Q_i)$. In other words, we have that
\[
 \sum_{i=1}^\infty \frac{v_i^2(t)}{Q_i i^2} = \sum_{i=1}^\infty Q_i h_i^2(t) < Ce^{Ct}.
\]
As the $Q_i$ decay exponentially, see \eqref{QLimit}, for any $T > 0$, and any $t \in [0,T]$ we have that
\[
|(\mathbb{L}V)_i| \leq \left|-a_iQ_1 v_i(t) + b_{i+1} v_{i+1}(t) \frac{i}{i+1}\right| \leq C i^{2+\alpha} Q_i,
\]
which goes to zero as $i \to \infty$. This implies that, uniformly for $t\in[0,T]$,
% for any fixed $t>0$
\[
\lim_{i \to \infty} |V_i(t) - V_i(0)| \leq \lim_{i \to \infty} \int_0^t |(\mathbb{L}V(s))_i| \,ds = 0.
\]
 As $\lim_{i \to \infty} V_i(0) = 1$ (since $\sum_{i=1}^\infty Q_i i h_i^0 = 1$), 
we then have that uniformly for $t\in[0,T]$, %, for any $t>0$, that 
\begin{equation}\label{eqn:v-limit}
\lim_{i\to \infty} V_i(t) = 1.
\end{equation}

{\bf Step 1: Minimum principle for $\mathbb{L}$.}
We claim that $\mathbb{L}$ satisfies a minimum principle in the sense that if
\begin{equation}\label{eqn:super-sol-def}
\partial_t W_i - (\mathbb{L} W)_i \geq 0, \quad \text{ for all } i \in \mathbb{N} ,\  t \in [0,T],
\end{equation}
then for any $N$
\[
\min_{i \in 1 \dots N, t \in [0,T]} W_i(t) \geq \min \left(\min_{i \in 1 \dots N} W_i(0), \min_{t \in [0,T]} W_N(t), \min_{t \in [0,T]} W_{1}(t)\right).
\]
To prove this, suppose that there exists some $j \in 2 \dots N-1$ and $\hat t \in (0,T]$ so that $W_{j}(\hat t) \leq W_i(t)$ for all $i \in 1 \dots N$ and $t \in [0,T]$. Then $\partial_t W_{j} (\hat t) \leq 0$, which in turn implies that $(\mathbb{L} W(\hat t))_{j} \leq 0$. However, the form of the differences in $\mathbb{L}W$ (see equation \eqref{eqn:integrated-form}), and the fact that $W_{j}(\hat t)$ is a minimizer, in turn implies that $W_{j -1}(\hat t) = W_{j +1}(\hat t) = W_{j}(\hat t)$. By repeating this argument, we find that $W_N(\hat t) = W_1(\hat t)= W_{j}(\hat t)$, which implies that the bound holds.

Exactly the same proof gives a maximum principle as well, in the sense that if 
\[
\partial_t W_i - (\mathbb{L} W)_i \leq 0, \quad \text{ for all } i\in\mathbb{N} ,\ t \in [0,T]
\]
then for any $N$
\[
\max_{i \in 1 \dots N, t \in [0,T]} W_i(t) \leq \max \left(\max_{i \in 1 \dots N} W_i(0), \max_{t \in [0,T]} W_N(t), \max_{t \in [0,T]} W_{1}(t)\right).
\]

{\bf Step 2: Supersolutions.} 
Define the function
\begin{equation}\label{eqn:def-w}
W^1(x,t) := \begin{cases}  \exp\left( \frac{x-A(N_1,2t)}{D}\right) &\text{ for } x < A(N_1,2t), \\ 1 &\text{ otherwise,} \end{cases} \\
\end{equation}
where $D>0$ is a constant that will be determined later.

We then claim that there exists some $N^*$ (independent of $N_1$) so that $W_i^1(t) := W^1(i,t)$ is a supersolution (i.e., \eqref{eqn:super-sol-def} is satisfied) for all $t>0$ satisfying $A(N_1,2t) > N^*$.
Clearly for $i > A(N_1,2t) + 1$, $(\mathbb{L} W^1)_i = 0$ and $\partial_t(W_i^1) = 0$, and hence \eqref{eqn:super-sol-def} is satisfied trivially. In the region $i < A(N_1,2t) - 1$ we compute:

\begin{align*}
\partial_t W_i^1 - (\mathbb{L} W^1)_i &= \frac{2(z_s-z)A(N_1,2t)^\alpha W_i^1}{D} + a_i Q_1 W_i^1 (1-e^{-D^{-1}}) - W_i^1\frac{b_{i+1} i}{i+1}(e^{D^{-1}} - 1) \\
&\geq  W_i^1 \left(\frac{2(z_s-z)A(N_1,2t)^\alpha}{D} + \frac{a_iQ_1 - \frac{b_{i+1}i}{i+1}}{D} -  \frac{C i^\alpha}{D^2}\right) \\
& = W_i^1 \left(\frac{2(z_s-z)A(N_1,2t)^\alpha}{D} + \frac{i^\alpha(z-z_s)}{D}+\frac{o(i^\alpha)}{D} - \frac{C i^\alpha}{D^2}\right),
\end{align*}
where on the last line we have used the assumption that $\fract{a_i}{i^\alpha} \to 1$, as well as \eqref{tildeQLimit}. We note that as long as $A(N_1,2t) > N^*$, with $N^*$ independent of $N_1,N_2$
(and $D$), the $o(i^\alpha)$ term will be dominated by $\fract{(z_s-z)A(N_1,2t)^\alpha}{2D}$ for all $i < A(N_1,2t)$. Hence for $D$ chosen large enough that $\frac{C}{D} < \frac{z_s - z}{2}$, and for $i < A(N_1,2t) - 1$  we will have that
\[
\partial_t W_i^1 - (\mathbb{L} W^1)_i \geq 0.
\]
It only remains to prove the boundary cases. If $A(N_1,2t) < i < A(N_1,2t)+1$ then
\[
\partial_t W_i^1 - (\mathbb{L} W^1)_i = - (\mathbb{L} W^1)_i = a_iQ_1(W^1_i - W^1_{i-1}) \geq 0.
\]
In the case $A(N_1,2t)-1 < i < A(N_1,2t)$ then
\begin{align*}
\partial_t W^1_i &- (\mathbb{L} W^1)_i \\
&= \frac{2(z_s-z)A(N_1,2t)^\alpha W^1_i}{D} + a_i Q_1 W^1_i (1-e^{-D^{-1}}) - \frac{b_{i+1} i}{i+1}(W^1_{i+1} - W^1_i) \\
&\geq \frac{2(z_s-z)A(N_1,2t,N_1)^\alpha W^1_i}{D} + a_i Q_1 W^1_i (1-e^{-D^{-1}}) - W^1_i \frac{b_{i+1} i}{i+1}(e^{D^{-1}} - 1) \\
& = W^1_i \left(\frac{2(z_s-z)A(N_1,2t)^\alpha}{D} + \frac{i^\alpha(z-z_s)}{D}+\frac{o(i^\alpha)}{D} - \frac{C i^\alpha}{D^2}\right) \geq 0.
\end{align*}
Hence $(W^1)$ is a supersolution in the sense that it satisfies \eqref{eqn:super-sol-def}. 

Similarly, if we define
\[
W^2(x,t) := \begin{cases} 1 &\text{ if } x < A(N_2,t/2), \\ 
\exp\left( \frac{A(N_2,t/2)-x}{D}\right) &\text{ otherwise}.\end{cases}
\]
we can use exactly the same type of estimates to show that that $W_i^2(t) := W^2(i,t)$ is also a supersolution in the sense that \eqref{eqn:super-sol-def} is satisfied, as long as $A(N_1,2t) > N^*$.

{\bf Step 3: Support bounds for $v$.}
We claim that $W^1_i(t) \geq V_i(t)$ for all $i \in \mathbb{N}$ and $t>0$ such that $A(N_1,2t) > N^*$. To prove this, we first note that because $A(N_1,0) = N_1$, we have that $W^1(0) \geq V(0)$. Furthermore, by \eqref{eqn:v-limit} and \eqref{eqn:def-w} we have that 
\[
\lim_{i \to \infty} W_i^1(t) -V_i(t) = 0,
\]
for all $t>0$ such that $A(N_1,2t) > N^*$. Since $W^1 - V$ is a supersolution, the minimum principle then implies that
\[
\min_{i \in \mathbb{N}, t \in [0,T]} W_i^1(t) - V_i(t) \geq \min \left(0,\min_{t \in [0,T]} W_1^1(t) - V_1(t) \right).
\]
Suppose for the sake of contradiction that 
\[
    \min_{t \in [0,T]} W_1^1(t) - V_1(t) = W_1^1(\hat t)-V_1(\hat t) < 0
\]
for some $\hat t \in (0,T]$. Clearly $\partial_t(W_1^1 - V_1)(\hat t) \leq 0$, and as $W^1 - V$ is a supersolution then $(\mathbb{L}(W^1-V))_1 \leq 0$. The definition of $\mathbb{L}$, along with the fact that $W^1-V$ is minimized at $j=1$, $t = \hat t$, then gives that
\begin{align*}
0 &\geq (\mathbb{L}(W^1-V))_1 = -a_iQ_1 (W_1^1 - V_1) + b_{i+1} ((W_2^1 - V_2) - (W_1^1 - V_1)) \frac{i}{i+1}\\
 &\geq -a_iQ_1 (W_1^1 - V_1) > 0,
\end{align*}
which is a contradiction. This then implies that
\[
\min_{i \in \mathbb{N}, t \in [0,T]} W_i^1(t) - V_i(t) \geq 0,
\]
which is the desired conclusion.

This readily implies that for any $\e$ there exists a $K^*$ (which depends only upon $D$ and $\e$) so that for all $t$ small enough that $A(N_1,2t) > N^*$ we have that
\[
\sum_{i=1}^{\lceil A(N_1,2t)-K^* \rceil} v_i(t) = V_{\lceil A(N_1,2t)-K^* \rceil} \leq W_{\lceil A(N_1,2t)-K^* \rceil}^1 < \frac{\e}{2}.
\]

On the other hand, $W^2(0) \geq 1-V(0)$. Hence by the minimum principle, $W^2(t) \geq 1-V(t)$ for all $t>0$ such that $A(N_1,2t) > N^*$, which thus implies, for all such $t$, that
\[
\sum_{i=\lceil A(N_2,t/2)+K^*\rceil}^\infty  v_i(t) < \frac{\e}{2}.
\]
Next, we observe that the coefficient of $v_i$ in $(\tfL v)_i$ is negative, whereas all of the other coefficients are positive, see \eqref{eqn:tilde-L}. This readily implies that if $v^0 \geq 0$ then $v_i(t) \geq 0$ for all $i \in \mathbb{N}, t \geq 0$. In turn, we may use the previous inequality to deduce that for all $t>0$ such that $A(N_1,2t) > N^*$,
\[
\|v - \chi v\|_{\ell^1} < \e \,.
\]

{\bf Step 4 Duhamel estimates:} By Duhamel's formula, one has that 
\begin{equation}\label{e:Duhamelv}
v(t) = e^{\tfL t} u^0 = e^{\bfL t} u^0 + \int_0^t e^{\bfL (t-s)} (\tfL -\bfL ) v(s) \ds.
\end{equation}

We seek to derive some simple bounds on the last term, which will provide the estimates we need. To that end, we estimate
\begin{align*}
    & \left\| \int_0^t e^{\bfL (t-s)} (\tfL -\bfL ) v(s) \ds
    \right\|_{\ell^1} \leq M \int_0^t \|(\tfL - \bfL )v(s)\|_{\ell^1} \ds \\
    &\hspace{1cm}\leq\  M \int_0^t \sum_{i=2}^\infty  \left|v_1i(-a_{i-1}Q_{i-1} + a_iQ_i) +\frac{b_i}{i}v_i - \frac{a_{i-1}Q_1}{i-1}v_{i-1} \right|\ds \\
    &\hspace{1cm} \quad    + M \int_0^t \left| -a_{1} Q_1 v_1 
+ \sum_{i=1}^\infty \left(\frac{b_{i+1}}{i+1}v_{i+1} - \frac{a_iQ_1}{i} v_i - a_iQ_i v_1 \right) \right| \ds,
\end{align*}
where we have used that $e^{\bfL t}$ has uniformly bounded operator norm in $\ell^1$ (since $e^{Lt}$ is uniformly bounded in  $\ell^1(Q_i i)$). We remark that by our assumption that $\fract{a_i}{i^\alpha} \to 1$, and the fact that the $Q_i$ decay exponentially, we have that
\[
\sum_{i=2}^\infty  \left|v_1i(-a_{i-1}Q_{i-1} + a_iQ_i) +\frac{b_i}{i}v_i - \frac{a_{i-1}Q_1}{i-1}v_{i-1} \right| \leq C\sum_{i=1}^\infty i^{\alpha-1} |v_i|.
\]
Then using the facts that $V \leq W^1$, $v_i(t) \geq 0$, and that $0< \alpha <1$, we can estimate
\begin{align*}
\int_0^t \sum_{i=1}^\infty i^{\alpha-1} |v_i|  
&\leq \int_0^t \sum_{i=1}^{\lceil A(N_1,2s)/2 \rceil } i^{\alpha-1} |v_i| +  \sum_{\lceil A(N_1,2s)/2 \rceil}^\infty i^{\alpha-1} |v_i| \ds \\
&\leq \int_0^t \exp \left(\frac{-A(N_1,2s)}{2D}\right) + \left( \frac{A(N_1,2s)}{2} \right)^{\alpha-1} \ds.
\end{align*}
%----------

By the change of variables $s\mapsto A(N_1,2s)$, using \eqref{eqn:A-sol} to write
$2(z_s-z) ds = -A^{-\alpha} dA$, and assuming $t\le T$ where $N_T:=A(N_1,2T)>1$, we find that
\begin{align*}
\int_0^t \sum_{i=1}^\infty i^{\alpha-1} |v_i|  
&\leq 
C \int_{N_T}^{N_1} \left(  A^{-\alpha} \exp\left(\frac{-A}{2D}\right)  + A\inv \right) \, dA
\\ &\leq 
    C\left( N_T^{1-\alpha}\exp\left(\frac{-N_T}{2D}\right)
    \left(\frac{N_1}{N_T}-1\right)+ \log\left(\frac{N_1}{N_T}\right) \right)\,.
\end{align*}
%----------
%
%By then using equation \eqref{eqn:A-sol}, along with the fact that $0< \alpha < 1$, and assuming that $t$ is small enough that $A(N_1,2s) > 1$, we can write 
%\[
%    \exp \left(\frac{-A(N_1,2s)}{2D}\right) \leq \exp \left(\frac{-A(N_1,2s)^{1-\alpha}}{2D}\right),
%    \]
%and thus we find that
%\begin{align*}
%\int_0^t \sum_{i=1}^\infty i^{\alpha-1} |v_i|  \leq &\frac{D}{(z_s-z)(1-\alpha)} \exp\left(\frac{-N_1^{1-\alpha} + 2(z_s-z)(1-\alpha)t}{2D}\right) \\
%&+\frac{1}{2(z_s-z)(1-\alpha)} \log \left( \frac{N_1^{1-\alpha}}{N_1^{1-\alpha} - 2(z_s-z) (1-\alpha) t}\right).
%\end{align*}
Similarly, we note that
\[
\left| -a_{1} Q_1 v_1 
+ \sum_{i=1}^\infty \frac{b_{i+1}}{i+1}v_{i+1} - \frac{a_iQ_1}{i} v_i - a_iQ_i v_1   \right|  \leq C\sum_{i=1}^\infty i^{\alpha-1} |v_i|,
\]
and hence we find that
\[
\left\| \int_0^t e^{\bfL (t-s)} (\tfL -\bfL ) v(s) \ds \right\|_{\ell^1} \leq C\left(   \left( \frac{N_1}{N_T} -1\right)+ \log\left(\frac{N_1}{N_T}\right) \right)\,,
    \]
where $C$ is independent of $v$. Hence, recalling \eqref{eqn:A-sol}, we find that for any $\e>0$ there exists a $\delta$ so that for $t < \delta N_1^{1-\alpha}$ we have that
\[
\left\| \int_0^t e^{\bfL (t-s)} (\tfL -\bfL ) v(s) \ds \right\|_{\ell^1} < \e .
\]
In turn, by Duhamel's formula \eqref{e:Duhamelv}, this immediately implies that
for all $t< \delta N_1^{1-\alpha}$,
\[
\|u-v\|_{\ell^1} \leq \e\,.
\]

{\bf Step 5: Cutoff estimates for $u$.}
Using the triangle inequality, and the support estimates on $v$ from Step 3, we find that
\begin{align*}
\|u-\chi u\|_{\ell^1} &\leq \|u-v\|_{\ell^1} + \|v-\chi  v\|_{\ell^1} 
%\\ &
+ \|\chi (v-u)\|_{\ell^1} 
%\\ &
\leq 3\e
\end{align*}
for any $t < \delta N_1^{1-\alpha}$.
This, in light of \eqref{eqn:mass-preserved}, finishes the proof of Theorem~\ref{thm:dynamics}.
\endproof
%\end{proof}

Finally, 
we prove Corollary \ref{cor:cutoff}, as a natural consequence of Theorem \ref{thm:dynamics}.

%\begin{proof}
\textit{Proof of Corollary~\ref{cor:cutoff}.} Let
\begin{align*}
u_i^1 &= \begin{cases} 2/N &\text{ for } i \in [N/4,N/2],\\ 0&\text{ otherwise}.  \end{cases}\\
u_i^2 &= \begin{cases} 2/N &\text{ for } i \in [3N/4,N],\\ 0&\text{ otherwise}.  \end{cases}
\end{align*}
Let $\chi^1$ and $\chi^2$ be the functions associated with $u^1$ and $u^2$ from Theorem \ref{thm:dynamics}. By the form of $A$, we know that for some $\hat \delta>0$ independent of $N$ we know that $\chi^1 \chi^2 = 0$ for all $t< \hat \delta N^{1-\alpha}$.

Let $u^0 = u^1 - u^2$. By the form of $A$, that is \eqref{eqn:A-sol}, along with Theorem \ref{thm:dynamics}, we know that
\begin{align*}
\|e^{\bfL t}u^0\|_{\ell^1} &\geq \sum_{i=1}^\infty \chi^1 |(e^{\bfL t}(u^1 - u^2))_i| + \sum_{i=1}^\infty \chi^2 |(e^{\bfL t}(u^1-u^2))_i| \\
 &\geq \sum_{i=1}^\infty \chi^1 |(e^{\bfL t}u^1)_i| - \chi^1 |(e^{\bfL t}u^2)_i| +\chi^2 |(e^{\bfL t}u^2)_i| - \chi^2 |(e^{\bfL t}u^1)_i| \\
 & \geq 1-\e.
\end{align*}
This proves the first part of the corollary.

For the second part, we recall (see section 3 in \cite{CanizoLods}) that $e^{Lt}$ generates a semigroup on the space
\[
Y_\eta := \left\{ h : \|h\|_\eta < \infty,\quad \|h\|_\eta := \sum_{i=1}^\infty Q_i e^{\eta i } h_i,\quad  \sum_{i=1}^\infty Q_i i h_i = 0 \right\}
\]
as long as $\eta$ sufficiently small, and that $e^{Lt}$ satisfies the bound
\[
\|e^{Lt} h^0\|_\eta \leq C e^{-\lambda t} \|h^0\|_\eta
\]
for some $\lambda > 0$, and for all $h \in Y_\eta$.

    We then note that any element of $X_1$ with support in $\{i < N\}$ will also be an element of $Y_\eta$, and will satisfy
\[
\|h^0\|_\eta = \sum_{i=1}^N Q_i e^{i \eta} |h_i^0| \leq e^{N \eta} \sum_{i=1}^\infty Q_i i |h_i^0| = e^{N\eta}.
\]
Hence we find that
\[
\|e^{Lt}h^0\|_1 \leq C\|e^{Lt} h^0 \|_\eta \leq C e^{-\lambda t} \|h^0\|_\eta \leq C e^{-\lambda t + N \eta}.
\]
The result then follows.
\endproof
%\end{proof}
% 
% \section*{Acknowledgements}
% This material is based upon work supported by the National
% Science Foundation under grants 
% DMS 1211161 and DMS 1515400,
% and partially supported by the Center for Nonlinear Analysis (CNA)
% under National Science Foundation PIRE Grant no.\ OISE-0967140,
% and the NSF Research Network Grant no.\ RNMS11-07444 (KI-Net).

\bibliography{BDquasilinearRefs}
\bibliographystyle{siam}

 \end{document}